\documentstyle[12pt,leqno]{amsart}

\newtheorem{lemma}{Lemma}
\newtheorem{corollary}{Corollary}  
\newtheorem{proposition}{Proposition}

\newcommand{\Spec}{\operatorname{Spec}}
\newcommand{\Sym}{\operatorname{Sym}}
\newcommand{\Aut}{\operatorname{Aut}}
\newcommand{\Hom}{\operatorname{Hom}}
\renewcommand{\Im}{\operatorname{Im}}
\renewcommand{\Re}{\operatorname{Re}}
\newcommand{\SL}{\operatorname{SL}}
\newcommand{\ad}{\operatorname{ad}}

\newcommand{\sheafhom}{\mathop{\text{\it $\cal H$om}}}
\newcommand{\ip}[2]{\left\langle{#1}\vphantom{#2}
   \right.\left|\,\vphantom{#1}{#2}\right\rangle}
\newcommand{\mufive}{\mbox{\boldmath ${\mu}_5$}}
\newcommand{\suchthat}{\ | \ }
\newcommand{\exterior}{\Lambda }
\newcommand{\prim}[1]{{#1}'}

\begin{document}

\title[Mirror symmetry and rational curves on quintic threefolds]{
Mirror symmetry and rational curves on quintic \\threefolds:
A guide for mathematicians}

\author{David R. Morrison}

\address{
Department of Mathematics \\
Duke University \\
Durham, NC\ \ 27706
}

\email{drm@@math.duke.edu}

\maketitle

\setlength{\unitlength}{1in}
\hfill
\begin{picture}(2,0)
\put(0,2.2){\makebox{\parbox{2in}{\flushright
           {\sc duk-m-91-01} \\ July, 1991}}}
\end{picture}

\begin{abstract}
We give a mathematical account of a recent string theory  calculation
which predicts the number of rational curves on the generic
quintic threefold.  Our account involves
the interpretation of Yukawa couplings in terms of variations of
Hodge structure,
a new $q$-expansion principle
for functions on the moduli space of Calabi-Yau manifolds,
and the ``mirror symmetry'' phenomenon recently observed by string
theorists.
\end{abstract}

\section*{Introduction}

There has been much recent excitement among mathematicians
about a calculation made by a group of string theorists
(P.~Candelas, X.~C.~de la Ossa, P.~S.~Green, and L.~Parkes \cite{pair})
which purports to give a count of the number of rational curves
 of fixed degree on a general quintic threefold.
The calculation mixes arguments from string theory with arguments
from mathematics, and is generally quite difficult to follow for
mathematicians.

I believe that I now understand the essential
mathematical content of that calculation.
It is my purpose in this note to explain my understanding
in terms familiar to algebraic geometers.
What Candelas et al.\ actually calculate is a $q$-expansion of a
certain function determined by the variation of Hodge structure
of some {\em other\/} family of threefolds with trivial canonical
bundle.  The ``mirror symmetry principle'' is then invoked to predict
that the Fourier coefficients in that expansion should be related
to the number of rational curves on a quintic threefold.

One mathematical surprise in this story is a new $q$-expansion principle
for functions on the moduli space of Calabi-Yau manifolds.
Near points on the boundary of moduli where the monodromy is ``maximally
unipotent,'' there turn out to be
 natural coordinates in which to make $q$-expansions of functions.
In this paper, we will discuss these
 $q$-expansions only
 in the case of one-dimensional moduli spaces; the general
case will be treated elsewhere.

By focusing on this $q$-expansion principle, we place the computation
of \cite{pair} in a mathematically natural framework.  Although there
remain certain dependencies on a choice of coordinates, the coordinates
used for calculation are canonically determined by the monodromy of
the periods, which is itself intrinsic.  On the other hand, we have
 removed some of the physical arguments which were
used in the original paper to help choose the coordinates appropriately.
The result may be that our presentation is less convincing to physicists.

\bigskip

The plan of the paper is as follows.  In section
\ref{Yukawa} 
we review variations of Hodge structure, and explain how to define
``Yukawa couplings'' in Hodge-theoretic terms.
(A discussion of Yukawa couplings along the same lines has
also been given by Cecotti \cite{cecotti1}, \cite{cecotti2}.)
In section
\ref{period} 
we study the asymptotic behavior of the periods near points with maximally
unipotent monodromy.  This is applied to find $q$-expansions
of the Yukawa couplings in section
\ref{qexp}. 
In section \ref{mirror} 
we attempt to describe  mirror symmetry in geometric terms.  In section
\ref{qm} 
we turn to the main example (the family of ``quintic-mirrors''),
and in section
\ref{moonshine} 
we explain how the mirror symmetry principle predicts from the
earlier calculations what the numbers of rational curves on quintic
threefolds should be.  Several technical portions of the paper have
been banished to appendices.

We work throughout with algebraic varieties over the complex numbers,
which we often identify with complex manifolds (or complex analytic
spaces).
If $X$ is
 a compact complex manifold and $p$, $q\ge 0$, we define
\[
H^{p,q}(X)=H^q(\Omega^p_X)=H^q(\exterior^p\Omega _X)
\]
 where $\Omega _X$ is the holomorphic cotangent bundle of $X$.
(This is slightly non-standard.)
We extend this definition to the case $p<0$, $q\ge 0$ by
\[
H^{p,q}(X)=H^q(\exterior^{-p}\Theta_X),
\]
where $\Theta_X$ is the holomorphic tangent bundle of $X$.
(This is very non-standard.)
The dimension
of $H^{p,q}(X)$ is denoted by $h^{p,q}(X)$, or simply  by $h^{p,q}$.

\section{Variations of Hodge structure arising from families of
Calabi-Yau manifolds}
\label{Yukawa}

Recall that a {\em Calabi-Yau manifold} is a compact K\"ahler manifold $X$
of complex dimension $n$ which has trivial canonical bundle,
such that the Hodge numbers $h^{k,0}$ vanish for $0<k<n$.
Thanks to a celebrated theorem of Yau \cite{yau}, every such manifold
admits Ricci-flat K\"ahler metrics.

Any non-zero section of the canonical bundle determines
isomorphisms
\[
H^i(\Theta_X) \cong H^{i}(\Omega^{n-1}_X).
\]
Thus, if $n>1$, then $X$ has no holomorphic vector fields.  Moreover, the
tangent space to moduli $H^1(\Theta_X)$ has dimension
$h^{-1,1}=h^{n-1,1}$ and
the natural obstruction
space for the moduli problem $H^2(\Theta_X)$ has dimension $h^{n-1,2}$,
which is generally non-zero for $n>2$.
However,
the theorem of Bogomolov \cite{bogomolov}, Tian \cite{tian}, and
Todorov \cite{todorov}
says that the moduli problem is in fact unobstructed, and the moduli space
is therefore smooth of dimension $h^{n-1,1}$.

We review some facts about Hodge structures and their variation.
Good general references for this
 are
Griffiths et al.\ \cite{transcendental}, and Schmid \cite{schmid}.
The $n^{\text{th}}$ cohomology group of $X$ has a {\em Hodge decomposition}
\[
H^n(X,\Bbb C) \cong \bigoplus\begin{Sb} p+q=n\\ p,q\ge 0 \end{Sb} H^{p,q}(X).
\]
(With our conventions, this follows from the Hodge theorem in de Rham
cohomology $H^n_{\text{DR}}(X) = \bigoplus H^{p,q}_{\bar\partial}(X)$
together with the Dolbeault isomorphism
$H^{p,q}_{\bar\partial}(X)\cong H^q(\Omega^p_X) = H^{p,q}(X)$.)
The Hodge decomposition can also be described by means of the {\em Hodge
filtration}
\[
F^p(X):=\bigoplus_{p\le p'\le n} H^{p',n-p'}(X);
\]
we then have $H^{p,n-p}(X)\cong F^p(X)/F^{p+1}(X)$.

The cup product on cohomology composed with evaluation on
the canonical orientation class of $X$
 determines a bilinear map
\[
\ip{\,}{\,}: H^n(X,\Bbb Z) \times H^n(X,\Bbb Z) \to H^{2n}(X,\Bbb Z)
\overset{\cong}{\to} {\Bbb Z},
\]
called a {\em polarization}.
There is an associated {\em adjoint map}
\[
\ad_{\ip{\,}{\,}}: H^n(X,\Bbb Z) \to
\Hom(H^n(X,\Bbb Z),\Bbb Z)
\]
defined by
$
\ad_{\ip{\,}{\,}}(x)(y)=\ip{x}{y}
$.
After tensoring with ${\Bbb C}$ and invoking the Hodge decomposition,
the adjoint map induces isomorphisms
\[
\ad_{\ip{\,}{\,}}: H^{p,n-p}(X) \overset{\cong}{\to}
(H^{n-p,p}(X))^*.
\]

We now recall a construction which first arose in the study of
infinitesimal variations of Hodge structure by Carlson, M.~Green, Griffiths,
and Harris \cite{CGGH}.
The cup product determines a natural map
\begin{equation}  \label{eqF}
H^1(\Theta_X)\to \bigoplus \Hom(H^{p,q}(X), H^{p-1,q+1}(X))
\end{equation}
called the {\em differential of the period map}.  Iterates of this map
are symmetric in their variables; the {\em $n^{\text{th}}$ iterate
of the differential} is the induced map
\[  
\Sym^n H^1(\Theta_X)\to\Hom(H^{n,0}(X), H^{0,n}(X)) .
\]  
Using the canonical isomorphism $\Hom(H^{n,0}(X), H^{0,n}(X))=
(H^{n,0}(X))^*\otimes H^{0,n}(X)$
and the isomorphism
$H^{0,n}(X)\cong(H^{n,0}(X))^*$ induced by the adjoint,
we get a map
\begin{equation}
\label{eqB}
\Sym^nH^1(\Theta_X)\to(H^{n,0}(X)^*)^{\otimes 2}.
\end{equation}
We call this the {\em unnormalized Yukawa
coupling\/}\footnote{This particular Yukawa coupling
is probably only interesting in physics
if $n=3$.  In dimension $n$, what is being computed here is the
``$n$-point Yukawa coupling.''}
 of $X$.

If we choose an element
of $H^{n,0}(X)^{\otimes 2}$ and evaluate the map (\ref{eqB}) on that element,
we get a map
$
\Sym^nH^1(\Theta_X)\to\Bbb C
$
called a {\em normalized Yukawa coupling}.  The ``normalization'' is the
choice of element of $H^{n,0}(X)^{\otimes 2}$.

\bigskip

We now analyze these constructions for a {\em family} of manifolds.
Suppose we are given a quasi-projective variety $C$ and a
smooth map $\pi: {\cal X}\to C$ whose fibers are
 Calabi-Yau manifolds.  Suppose also that this family
 can be completed to a family of varieties
$\bar\pi: \overline{\cal X}\to \overline C$, where $\overline C$ is a
projective compactification of $C$.  The fibers of $\bar\pi$ may degenerate
over the boundary $B:=\overline C - C$.  We assume that $B$ is a divisor
with normal crossings on $\overline C$.

For any point $P\in C$ we denote the fiber $\pi^{-1}(P)$ of $\pi$
by $X_P$.  The Kodaira-Spencer map
$
\rho:\Theta_{C,P}\to H^1(\Theta_{X_P})
$
maps the tangent space to $C$ at $P$ to the tangent space to moduli at
the point $[X_P]$.  If $C$ is actually a moduli space for the fibers
of ${\cal X}$, the map $\rho$ will be an isomorphism.

The cohomology of the fibers of the map $\pi$ with coefficients in
${\Bbb Z}$ and ${\Bbb C}$ fit together into local
systems $R^n\pi_* {\Bbb Z}$ and $R^n\pi_* {\Bbb C}$ on $C$.
The Hodge filtration becomes
 a filtration of the vector bundle
${\cal F}^0 := R^n\pi_* {\Bbb C} \otimes {\cal O}_C$ by
holomorphic subbundles:
\[
R^n\pi_* {\Bbb C} \otimes {\cal O}_C
= {\cal F}^0 \supset {\cal F}^1 \supset \cdots \supset
{\cal F}^{n-1} \supset
{\cal F}^n \supset  (0).
\]
The vector bundle ${\cal F}^0$ has
a natural flat connection $\nabla: {\cal F}^0\to {\cal F}^0\otimes \Omega_C$
called the {\em Gau\ss-Manin connection},
whose horizontal sections determine the local system
$R^n\pi_* {\Bbb C}$.
The {\em Griffiths transversality property} says that
$
\nabla ( {\cal F}^p )\subset {\cal F}^{p-1} \otimes \Omega _C
$.

There is a natural extension of this setup over the boundary $B$,
which involves the sheaf $\Omega_{\overline C}(\log B)$ of
{\em logarithmic differentials}.  (That sheaf is locally generated
by $\Omega_{\overline C}$ and all elements of the form $df/f$,
where $f=0$ is a local equation of a local component of $B$.)  Although
the local system $R^n\pi_* {\Bbb C}$
cannot in general be extended across $B$ in a single-valued way, the
Hodge bundles  ${\cal F}^p$ {\em do} have natural extensions to bundles
$\overline{\cal F}^p$ on $\overline C$.  And the Gau\ss-Manin
connection $\nabla$ extends to a connection
$
\overline\nabla : \overline{\cal F}^0 \to \overline{\cal F}^0 \otimes
\Omega_{\overline C}(\log B)
$
which satisfies
$
\overline
\nabla ( \overline{\cal F}^p )\subset \overline{\cal F}^{p-1}
\otimes \Omega _{\overline C}(\log B)
$.
This restriction on the types of poles which $\overline\nabla$
may have along $B$ is equivalent to a requirement that
the connection $\nabla$ have
 ``regular singular points.''

The extended Gau\ss-Manin connection $\overline\nabla$ gives
rise to an ${\cal O}_{\overline C}$-linear map on the associated gradeds
\begin{equation}  \label{eqD}
\widetilde\nabla: \overline{\cal F}^p/\overline{\cal F}^{p+1} \to
(\overline{\cal F}^{p-1}/\overline{\cal F}^{p}) \otimes
\Omega _{\overline C}(\log B).
\end{equation}
To make contact with the $n^{\text{th}}$ iterate of the differential
and the Yukawa coupling, we introduce the sheaf
$\Theta_{\overline C}({-\log B})$  of {\em vector fields with logarithmic
zeros},
which is the dual of $\Omega _{\overline C}(\log B)$.
The map (\ref{eqD}) then induces the bundle version of (\ref{eqF})
\[  
\Theta_{\overline C}({-\log B})\to \bigoplus \sheafhom (
\overline{\cal F}^p/\overline{\cal F}^{p+1} ,
\overline{\cal F}^{p-1}/\overline{\cal F}^{p}) .
\]  
When this is iterated $n$ times, it produces a map
\begin{equation}
\label{eqC}
\Sym^n(\Theta_{\overline C}({-\log B}))\to
\sheafhom (\overline{\cal F}^n ,
\overline{\cal F}^{0}/\overline{\cal F}^{1}).
\end{equation}

The polarizations fit together into a
bilinear map of local systems
\[
\ip{\,}{\,}: R^n\pi_*{\Bbb Z} \times R^n\pi_*{\Bbb Z} \to R^{2n}\pi_*{\Bbb Z}
\overset{\cong}{\to} {\Bbb Z}
\]
whose adjoint map induces an isomorphism
$
\ad_{\ip{\,}{\,}}: ({\cal F}^{0}/{\cal F}^{1}) \to
({\cal F}^{n})^*
$.
This extends to a map of bundles
\begin{equation} \label{eqG}
\ad_{\ip{\,}{\,}}: (\overline{\cal F}^{0}/\overline{\cal F}^{1}) \to
(\overline{\cal F}^{n})^*.
\end{equation}
Using the canonical isomorphism
$
\sheafhom (\overline{\cal F}^n ,
\overline{\cal F}^{0}/\overline{\cal F}^{1}) =
(\overline{\cal F}^n)^* \otimes
(\overline{\cal F}^{0}/\overline{\cal F}^{1})
$
and composing the map (\ref{eqG}) with the map (\ref{eqC}),
we get
the {\em Yukawa map}
\[
\kappa: \Sym^n(\Theta_{\overline C}({-\log B})) \to
((\overline{\cal F}^{n})^*)^{\otimes 2}.
\]
If we also specify a section of
$(\overline{\cal F}^{n})^{\otimes 2}$,
we get a {\em normalized Yukawa map}
\[
\kappa^{\text{norm}}: \Sym^n(\Theta_{\overline C}({-\log B})) \to
{\cal O}_{\overline C}.
\]

Suppose that
 $C$ is actually the moduli space for the fibers of ${\cal X}$, so that
$\rho$ is an isomorphism.  If we compose $\rho^{-1}$
with a normalized Yukawa map $\kappa^{\text{norm}}$ we get
\[
\Sym^nH^1(\Theta_{X_P}) \overset{\rho^{-1}}{\longrightarrow}
\Sym^n(\Theta_{C,P}) \overset{\kappa^{\text{norm}}}{\longrightarrow}
{\cal O}_{C,P} = \Bbb C.
\]
In this way, we exactly recover the corresponding
normalized Yukawa coupling.

Candelas et al.\ \cite{pair}
typically compute the Yukawa coupling in local coordinates
(away from the boundary) as follows.  Suppose that $\dim C = 1$,
and that $\psi$ is a local coordinate defined in an open set $U\subset C$.
There is an induced vector field $d/d\psi$, which is
 a local section  of $\Theta_U$.
Choose a section\footnote{To avoid confusion with the cotangent bundle,
we denote this section by $\omega$ rather than $\Omega$.  However,
in appendix C below, we will revert to the notation $\Omega$ used in
\cite{pair}.}
 $\omega$ of ${\cal F}^n$ over $U$, and define
\[
\kappa_{\psi\dots\psi}
=\kappa(\frac{d}{d\psi},\dots,\frac{d}{d\psi})\cdot\omega^2.
\]
(The number of $\psi$'s in the subscript is $n$.)  This is a holomorphic
function on $U$.  If we alter $\omega$
by the gauge transformation $\omega \mapsto f\omega$,
then the Yukawa coupling
transforms as $\kappa_{\psi\dots\psi} \mapsto f^2\kappa_{\psi\dots\psi}$.
``Normalizing the Yukawa map'' is the same thing as ``fixing the gauge.''

Our primary goal will be to compute the asymptotic behavior of the
Yukawa map $\kappa$ in a neighborhood of the boundary $B$.

\section{The asymptotic behavior of the periods}
\label{period}

For simplicity of exposition,
we now specialize to the case in which $C$ is a curve.
Let $P\in B$ be a boundary point, and let $T_P$ be the monodromy of the local
system $R^n\pi_*{\Bbb Z}$ around $P$.
We regard $T_P$ as an element of $\Aut H^n(X_{P'},\Bbb Z)$, where
$P'$ is a point near $P$; $T_P$ is determined by analytic continuation
along a path which goes once
around $P$ in the counterclockwise direction.  By the monodromy theorem
\cite{monodromy}, $T_P$ is quasi-unipotent, which means that some power
$T_P^k$ is unipotent.  Moreover, the index of unipotency is bounded:
we have $(T_P^k-I)^{n+1}=0$.

We say that {\em $P$ is a point at which the monodromy is maximally
unipotent} if the monodromy $T_P$ is unipotent, and if $(T_P-I)^n\ne 0$.
(Thus, the index of nilpotency of $T_P-I$ is maximal.)
Since $T_P-I$ is nilpotent, we can define the logarithm of the monodromy
$N=\log(T_P)\in\Aut H^n(X_{P'},\Bbb Q)$
by a finite power series
\[
\log(T_P)=(T_P-I)-\frac{(T_P-I)^2}2 +\dots +(-1)^{n+1}\frac{(T_P-I)^n}n.
\]
(Rational coefficients are needed in cohomology since
rational numbers appear in the power series.)
$N$ is also a nilpotent matrix, with the
same index of nilpotency as $T_P-I$.

\begin{lemma}
\label{lemma1}
Let $\pi: {\cal X}\to C$ be a one-parameter family of varieties with
$h^{n,0}=1$.
Let $P\in B=\overline C-C$ be a boundary point at which the monodromy
on $R^n\pi_*{\Bbb Z}$ is maximally unipotent
and let $N$ be the logarithm of the monodromy.
Then the image of $N^{n}$ is a $\Bbb Q$-vector space of dimension one,
and the image of $N^{n-1}$ is a $\Bbb Q$-vector space of dimension two.
\end{lemma}

We defer the proof of this lemma to appendix A.

\bigskip

We say that a basis $g_0$, $g_1$ of $(\Im N^{n-1})\otimes \Bbb C \subset
H^n(X_{P'},\Bbb C)$
is an
{\em adapted basis} if $g_0$ spans $(\Im N^n)\otimes \Bbb C$.
(We have extended scalars to $\Bbb C$ since certain computational
procedures lead more naturally to complex coefficients.)
If $g_0$, $g_1$ is an adapted basis for $(\Im N^{n-1})\otimes \Bbb C$, then
by Poincar\'e duality, there are homology classes
$\gamma_0, \gamma_1\in H_n(X_{P'},{\Bbb C})$
such that $\ip{g_j}{\alpha }=\int_{\gamma_j}\alpha $ for any
$\alpha \in H^n(X_{P'},{\Bbb C})$.
Here we denote the evaluation of cohomology classes on homology classes
by using an integral sign,
since that evaluation is often accomplished by integration.

\begin{proposition}
Let $\gamma_0, \gamma_1$ be the homology classes determined by an
adapted basis $g_0$, $g_1$ of $(\Im N^{n-1})\otimes \Bbb C$.
Define a constant $m$ by
$Ng_1=mg_0$.
Let $\overline U$ be a small neighborhood of $P$, and let $z$ be a
coordinate on $\overline U$ centered at $P$.  Let $\omega$ be a non-zero
section of $\overline{\cal F}^n$ over $\overline U$.
Then
\begin{enumerate}
\item
$\int_{\gamma_0}\omega$ extends to a single-valued function on
$\overline U$.

\item
$\int_{\gamma_1}\omega$ is not single-valued.  However, we can
write
\[
\frac{\frac1m\int_{\gamma_1}\omega}{\int_{\gamma_0}\omega}
= \frac{\log z}{2\pi i}  + \text{ single valued function.}
\]

\end{enumerate}
\end{proposition}

\begin{pf}
Any $g\in H^n(X_{P'},{\Bbb C})$ can be extended to a section $g(z)$ of
the local system over $U=\overline U - P$, which may be multi-valued.
But by the nilpotent orbit theorem  \cite{schmid},
$\exp(-\frac{\log z}{2\pi i}N)g(z)$ extends to a single-valued section.

Since $\omega$ is single-valued,
\[
\ip{\exp(-\frac{\log z}{2\pi i}N)g(z)}{\omega}
\]
will also be single-valued.  Now $g_j\in (\Im N^{n-1})\otimes \Bbb C$
implies that $N^2g_j=0$ for $j=1$, $2$.
The series needed for $\exp$ in this case is thus rather simple:
\begin{eqnarray*}
\exp(-\frac{\log z}{2\pi i}N)g_j(z) &=& (I - \frac{\log z}{2\pi i}N)g_j(z)\\
 &=& g_j(z) - \frac{\log z}{2\pi i}Ng_j(z).
\end{eqnarray*}

We conclude that
\[
\int_{\gamma_0}\omega = \ip{g_0(z)}{\omega}
\]
and
\[
\int_{\gamma_1}\omega - m\,\frac{\log z}{2\pi i}\int_{\gamma_0}\omega
=\ip{g_1(z) - \frac{\log z}{2\pi i}\,mg_0(z)}{\omega}
\]
are single-valued functions.
\end{pf}

\begin{corollary}
Let $\gamma_0, \gamma_1$ be the homology classes determined by an
adapted basis $g_0$, $g_1$ of $(\Im N^{n-1})\otimes \Bbb C$,
as in the proposition.  The function
\[
t:= \frac{\frac1m\int_{\gamma_1}\omega}{\int_{\gamma_0}\omega}
\]
gives a natural parameter
 on the universal cover $\widetilde U$
of  $U$ called a {\em quasi-canonical parameter}, and
\[ q:=e^{2\pi it} \]
gives a natural coordinate on $\overline U$ called a {\em quasi-canonical
coordinate}.  These are independent of
the choice of $\omega$.
We have
\[
\frac d{dt}=2\pi i\,q\,\frac d{dq},
\]
either of which serves as a local generator of the sheaf
$\Theta_{\overline C}({-\log B})$.

Moreover, under a change of adapted basis $(g_0,g_1)\mapsto (ag_0,bg_0+cg_1)$,
we have
\begin{eqnarray*}
m & \mapsto & \frac ca\,m, \\
t & \mapsto & t+\frac b{mc}, \text{\ and} \\
q & \mapsto & e^{2\pi ib/mc}q.
\end{eqnarray*}
Therefore, $t$ is uniquely determined up to an additive constant,
and $q$ is uniquely determined up to a multiplicative constant.
\end{corollary}

We can normalize further if we take the integral structure into
account.  We call $g_0$, $g_1$ a {\em good integral basis} of $\Im N^{n-1}$
if $g_0$ is a generator of
$\Im N^n \cap H^n(X_{P'},{\Bbb Z})$, and $g_1$ is an indivisible element
of $H^n(X_{P'},{\Bbb Z})$ which can be written as
$g_1 = \frac1\lambda N^{n-1}g$ for some $\lambda >0$ and some
 $g\in H^n(X_{P'},{\Bbb Z})$  such that $\ip{g_0}{g}=1$.
Notice that a good integral basis is an adapted basis.

The next lemma, which is based on some work of Friedman and Scattone
\cite{fs}, will be proved in appendix A.

\begin{lemma}
\label{lemma2}
Good integral bases exist.  If $g_0$, $g_1$ and $g_0'$, $g_1'$
are good integral bases, then
\begin{eqnarray*}
g_1' & = & k\, g_0 + (-1)^\ell g_1 \\
g_0' & = & (-1)^\ell g_0,
\end{eqnarray*}
for some integers $k$ and $\ell$.
\end{lemma}

Since $(T-I)^2 =0$ on $\Im N^{n-1}$, we have the simple formula $N=T-I$
on that space.  In particular, when restricted to
$\Im N^{n-1}$, the map $N$ is defined over the integers.  Thus,
if $g_0$, $g_1$ is a good integral basis and we write $Ng_1=mg_0$,
then $m$ is an integer.  Note that $m$ is independent of the choice
of good integral basis.

\begin{corollary}
Let $g_0$, $g_1$ be a good integral basis, and define
an integer $m$ by $Ng_1=mg_0$.
Then the quasi-canonical coordinate $q$ formed
from this basis  is unique up to multiplication by an $|m|^{\text{th}}$
root of unity.
\end{corollary}

We call $q$  a {\em canonical coordinate\/} and $t$ a {\em canonical
parameter\/} under these circumstances.
These are actually unique if $|m|=1$; in this
case, we say that the monodromy is {\em small}.

\section{The $q$-expansion of the Yukawa coupling}
\label{qexp}

The first example of the construction of the previous section is
furnished by the classical theory of periods of elliptic curves.
Let $\pi:{\cal X}\to U$ be a family of smooth elliptic curves over
a punctured disk $U$ which can be completed to a family
$\bar\pi: \overline{\cal X}\to \overline U$ with a singular
fiber over the boundary point $P=\overline U - U$.  The point $P$
is called a {\em cusp}.

Let $P'\in U$,
and suppose there is a symplectic
basis $\gamma_0$, $\gamma_1$ of the first homology
group $H_1(X_{P'},\Bbb Z)$ such that the monodromy $T_P$ acts as
\begin{eqnarray*}
T_P(\gamma_0) & = & \gamma_0 \\
T_P(\gamma_1) & = & \gamma_0 + \gamma_1.
\end{eqnarray*}
(The basis is {\em symplectic\/} if $\gamma_0\cap\gamma_1=1$.)
This easily implies that $P$ is a maximally unipotent boundary point,
that $\gamma_0$, $\gamma_1$ is the homology basis dual to a good
integral basis,
and that $m=1$.

For a fixed holomorphic one-form $\omega$ on $X_{P'}$, the numbers
$(\int_{\gamma_0}\omega, \int_{\gamma_1}\omega)$
were classically known
as the {\em periods\/} of the elliptic curve $X_{P'}$.  By varying the
one-form, the periods can be normalized to take the form
$(1,\tau)$.  The invariant way to formulate   this is to define
\[
\tau=\frac{\int_{\gamma_1}\omega}{\int_{\gamma_0}\omega}.
\]
This function $\tau$ can be regarded as a map from the universal cover
$\widetilde U$ of $U$ to the upper half-plane $\Bbb H$.  (The image lies
in the {\em upper\/} half-plane since the basis is symplectic.)

The monodromy transformation $T_P$ induces the map
\begin{equation} \label{fourier}
\tau\mapsto\tau+1.
\end{equation}
Thus, functions $f$ defined on $U$ pull back to functions $\widetilde f$
on $\widetilde U$
which are invariant under (\ref{fourier}).  It follows that any such
function has a {\em Fourier series\/}
\[
\widetilde f(\tau) = \sum_{n\in\Bbb Z}a_ne^{2\pi in\tau}.
\]
If expressed in terms of the natural coordinate $q=e^{2\pi i\tau}$ on
$U$, this is called a {\em $q$-expansion}, and it takes the form
\[
f(q)=\sum_{n\in\Bbb Z}a_nq^n.
\]
If $f$ has a holomorphic extension across the cusp $P$, the only terms
appearing in this sum are those with $n\ge 0$.

\bigskip

What we have shown in section \ref{period} is that this classical
construction
generalizes to functions defined near a maximally unipotent boundary point
$P$
of a Calabi-Yau moduli space (at least when that space has dimension one).
Fix a good integral basis, which determines a canonical coordinate $q$
and a canonical parameter $t$.
The monodromy transformation $T_P$ acts on $t$
by $t\mapsto t+1$.  Therefore, any function $f$ defined near $P$ which
is holomorphic at $P$
will have a $q$-expansion
\[
f(q)=\sum_{n= 0}^{\infty}a_nq^n,
\]
which can also be regarded as a Fourier series
\[
\widetilde f(t) = \sum_{n= 0}^{\infty}a_ne^{2\pi int}
\]
in $t$.
These expressions are unique if $|m|=1$, i.e., if the monodromy is small.

In order to obtain a $q$-expansion of the Yukawa coupling, we must
normalize that coupling.  But there is a natural choice of normalization
determined by a good integral basis.  To see this, note that
any good integral basis $g_0$, $g_1$  determines a section
$
(\int_{\gamma_0})^{-1} \in H^0(\overline U,(\overline{\cal F}^n))
$
by
\[
(\int_{\gamma_0})^{-1} := \frac{\omega}{\int_{\gamma_0}\omega}
\]
for any non-zero $\omega\in H^0(\overline U,\overline{\cal F}^n)$.
By lemma \ref{lemma2},
a change in good integral basis may change the sign of
$(\int_{\gamma_0})^{-1}$, but the induced section
$
(\int_{\gamma_0})^{-2}\in H^0(\overline U,(\overline{\cal F}^n)^{\otimes 2})
$
is independent of the choice of good integral basis.

We thus have a very natural normalization for the Yukawa map in $\overline U$.
We also have a natural parameter $t$ with which to compute, such that
$d/dt$ is a generator of $\Theta_{\overline U}({-\log B})$.
So we can define the {\em mathematically normalized
 Yukawa coupling\/} $\kappa_{t\dots t}$
by the formula
\[
\kappa_{t\dots t}
=\kappa(\frac d{dt},\dots,\frac d{dt})\cdot(\int_{\gamma_0})^{-2}.
\]

This mathematically normalized Yukawa coupling
$\kappa_{t\dots t}$ is an intrinsically defined function on a neighborhood
of the boundary.
(It is canonically determined by our
choice of maximally unipotent boundary point; however, it could
conceivably change if the boundary point changes.)
The function $\kappa_{t\dots t}$ therefore has a $q$-expansion
\begin{equation}\label{qexpq}
\kappa_{t\dots t} = a_0+a_1q+a_2q^2+\cdots,
\end{equation}
 which can also be regarded as a Fourier expansion in the
parameter $t$:
\begin{equation}\label{qexpt}
\kappa_{t\dots t} = a_0+a_1e^{2\pi it}+a_2e^{4\pi it}+\cdots.
\end{equation}
These expressions are unique if the monodromy is small.
\section{Mirror symmetry}
\label{mirror}

In this section I will attempt to outline the mirror symmetry
principle in mathematical terms, and describe some of the
mathematical evidence for it.  I apologize to
physicists
for my
misrepresentations of their ideas,
and I apologize to
mathematicians for the vagueness
of my explanations.

Gepner \cite{gepner1} has conjectured that there is a one-to-one
 correspondence between
$N=2$
superconformal field theories
with $c=3n$, and Calabi-Yau manifolds $X$ of dimension $n$
equipped with some
``extra structure'' $S$.
(This correspondence can be realized concretely in a number of
important cases using work of Greene, Vafa and Warner \cite{gvw},
Martinec \cite{martinec1}, \cite{martinec2}, and others.)
A precise geometric description of the
extra structure $S$ has not yet been given.  It appears to involve
specifying a class in $U/\Gamma$, where $U\subset H^{1,1}(X)$
is some open set, and $\Gamma$ is some group of automorphisms of
$U$.
What {\em is} clear about
this extra structure is how to perturb it:  first-order deformations
of $S$ correspond to elements of $H^{1,1}(X)$.

An instructive example is the case in which $X$ is
an  elliptic
curve.  In that case, as shown in \cite{dvv} and \cite{al},
one takes $U\subset H^{1,1}(X)\cong \Bbb C$
to be the upper half-plane, and $\Gamma=\SL(2,\Bbb Z)$.
Thus, the extra structure
$S$ represents a point in the $j$-line, or equivalently,
a choice of a {\em second} elliptic curve.

We specialize now to the case of dimension $n=3$. The space of first-order
deformations of the superconformal field theory can be decomposed
as\footnote{It has become common in the physics literature to use
$H^{2,1}(X)$ in place of $H^1(\Theta_X)$, largely because of the
success of Candelas \cite{candelas} and others in computing Yukawa
couplings on $H^{2,1}$.  In order to get the correct answer in families,
however, we must return to the original analysis of Strominger and
Witten \cite{str-wit} and work with Yukawa couplings on $H^1(\Theta_X)$.
The point is that while $H^1(\Theta_X)$ and $H^{2,1}(X)$ are isomorphic
for a Calabi-Yau threefold, they are not {\em canonically\/} isomorphic.
This affects the bundles over the moduli space to which they belong.}
 $H^1(\Theta_X) \oplus H^1(\Omega _X)$,
with
$H^1(\Theta_X)=H^{-1,1}(X)$ corresponding to first-order
deformations of the complex structure on $X$, and
$H^1(\Omega _X)=H^{1,1}(X)$ corresponding to first-order deformations of
the extra structure $S$.
These first-order deformations are called {\em marginal operators\/}
in the physics literature.

Specifying a superconformal field theory of this type also determines
 cubic forms $\Sym^3H^{-1,1}(X)\to\Bbb C$ and
$\Sym^3H^{1,1}(X)\to\Bbb C$.  The cubic form on $H^{-1,1}$ is the
Yukawa coupling described in section \ref{Yukawa}, normalized in
a way specified by the physical theory.  From a
mathematical point of view, this is determined by the variation of Hodge
structure plus the choice of normalization.
This cubic form depends on the complex structure of $X$, but should be
 independent of
the ``extra structure'' $S$.

The cubic form on $H^{1,1}$
lacks a precise geometric description at present.
By work of Dine, Seiberg, Wen, and Witten \cite{dsww2}
and Distler and Greene \cite{distler-greene},
it is known to have an expression of the form
\begin{equation}
\label{asymp}
 \sum_{k=0}^\infty  \sigma_k \, e^{-kR} ,
\end{equation}
where $R$ is a complex parameter which depends on the extra structure $S$.
The real part of $R$ is related to the ``radius'' in the physical theory
in such a way that $\Re R\to\infty$ is the ``large radius limit.''
The leading coefficient $\sigma_0$ is the natural topological product
$\Sym^3H^{1,1}(X)\to\Bbb C$.
(In other words, the cubic form on $H^{1,1}$ approaches the topological
product in the large radius limit.)
The higher coefficients $\sigma_k$ are
supposed to be related in some
well-defined way to the numbers of rational curves of various degrees
on the generic deformation of $X$ (assuming those numbers are finite).
One of the important unsolved problems in the theory is to determine this
relationship precisely.

As was first noticed by Dixon \cite[p.\ 118]{dixon}, and later developed
by Lerche, Vafa, and Warner \cite{lvw} and others, the identification
of one piece of the superconformal field theory with $H^{1,1}(X)$
and the other piece with $H^{-1,1}(X)\cong H^{2,1}(X)$ involves an arbitrary
choice, and the theory is also consistent with making the opposite choice.
Moreover, as we will describe below, there are examples in which the
Gepner correspondence can be realized for both choices.
But except in the very rare circumstance that
 the Hodge numbers $h^{1,1}$ and $h^{2,1}=\dim H^{-1,1}$
 coincide, changing the choice necessarily involves changing the
Calabi-Yau threefold $X$.  The new threefold $\prim X$ will
have a completely different topology from the old:  in fact, the
Hodge diamond is rotated by $90^\circ$ when passing from one
to the other.

This leads to a mathematical version of the mirror symmetry conjecture:
To each pair $(X,S)$ consisting
of a Calabi-Yau
threefold $X$ together with some extra structure $S$ there
should be associated a ``mirror pair''
$(\prim X,\prim S)$ which comes equipped with natural isomorphisms
$H^{-1,1}(X)\overset{\cong}{\to}H^{1,1}(\prim{X})$
and
$H^{1,1}(X)\overset{\cong}{\to}H^{-1,1}(\prim{X})$
that are compatible with the cubic forms.

Even in this rather imprecise\footnote{Among the things not properly
defined from a mathematical viewpoint, we must include
the normalization of the Yukawa coupling, the complex parameter $R$
(which depends on the ``extra
structure'' $S$), and the higher coefficients
$\sigma_k$.}
form, the conjecture as stated
is easily refuted:
There exist rigid Calabi-Yau threefolds, which have $h^{2,1}=0$
(see Schoen \cite{schoen}
for an example).  Any mirror of such a threefold
would have $h^{1,1}=0$, and so could not be
K\"ahler.
A potentially correct version of the conjecture,
even less precise, begins:  ``To most pairs
$(X,S)$, including almost all of interest in physics, there should be
associated\dots''.

It is tempting to speculate that the theory should
be extended to non-K\"ahler threefolds as in Reid's fantasy
\cite{nevertheless}, which might rescue the conjecture in its
original form.
Alternatively, Aspinwall and L\"utken \cite{al2} suggest that the
Gepner correspondence (and hence the mathematical version of
mirror symmetry) should only hold in the large radius
limit.  Since no ``limits'' can be taken in the rigid case,
a mathematical mirror construction would not be expected there.

To be presented with a conjecture which has been only vaguely formulated
is unsettling to many mathematicians.
Nevertheless, the
mirror symmetry phenomenon appears to be quite widespread,
so it seems important to make further efforts to find a precise
formulation.  In fact, there are at least
four major pieces of mathematically significant evidence for mirror symmetry.
\begin{enumerate}
\item[(i)]
Greene and Plesser \cite{greene-plesser} have studied a case in which
there are very solid physics arguments which tie the pair $(X,S)$
to the corresponding superconformal field theory (as predicted by Gepner).
The Calabi-Yau threefolds in question are desingularizations of
 quotients of Fermat-type
weighted hypersurfaces by certain finite groups (including the trivial
group).  For each pair $(X,S)$ of
this type, Greene and Plesser were able to find the corresponding
mirror pair $(\prim X,\prim S)$
by analyzing the associated superconformal field theories.
It turns out that the pairs are related by taking quotients:
$\prim X$ is a desingularization of $X/G$ for some symmetry group $G$.
By deformation arguments, the mirror symmetry phenomenon persists in
neighborhoods of $(X,S)$ and $(\prim X,\prim X)$.
 Roan \cite{roanpf} subsequently gave a
direct mathematical proof that the predicted isomorphisms between
$H^{-1,1}$ and $H^{1,1}$ groups exist in this situation.

\item[(ii)]
Candelas, Lynker, and Schimmrigk \cite{cls} have computed the Hodge
numbers for a large class of Calabi-Yau threefolds which are
desingularizations of hypersurfaces in weighted
projective spaces.  They put some extra constraints on the form of the
equation, and found about 6000 types of threefolds satisfying their
conditions.
The set of pairs $(h^{1,1},h^{2,1})$ obtained from these examples
is very nearly (but not precisely)
symmetric with respect to the interchange $h^{1,1}{\leftrightarrow}h^{2,1}$.
Since there is no {\em a priori} reason that the mirror of a desingularized
weighted hypersurface should again be a desingularized weighted
hypersurface, this is consistent with the conjecture and supports it
quite strongly.

\item[(iii)]
Aspinwall, L\"utken, and Ross \cite{alr} (see also \cite{al})
have carefully studied
a particular mirror pair $(X,S)$, $(\prim X,\prim S)$.
They put $X$ in a family ${\cal X} = \{X_t \}$ which has a degenerate limit
as $t$ approaches $0$.  Some heuristics were used in choosing the
family $\cal X$, in an attempt to ensure that the limit as $t\to 0$
would correspond to the ``large radius limit'' for the mirror
($\prim X,\prim S)$.  Aspinwall et al.\ then computed
the limiting behavior of
the cubic form on $H^{-1,1}(X_t)$, and showed that it coincides
with the topological product $\prim \sigma_0$ on
$H^{1,1}(\prim X)$,
as predicted by the conjecture.  (Actually, there is a normalization
factor which was not computed, but the agreement is exact up to
this normalization.)

\item[(iv)]
The work of Candelas,
de la Ossa, P.~Green, and Parkes \cite{pair} being described in this paper
goes further, and computes the other coefficients in an asymptotic
expansion.  This will be explained in more detail in the next two sections.
\end{enumerate}

\section{The quintic-mirror family}
\label{qm}

We now describe a certain one-parameter family of Calabi-Yau threefolds
constructed by Greene and Plesser \cite{greene-plesser}, as amplified
by Candelas et al.\ \cite{pair}.
Begin with the family of quintic threefolds
${\cal Q}_\psi = \{\vec{x}\in {\Bbb P}^4\suchthat p_\psi(\vec{x})=0\}$
defined by the polynomial
\[
p_\psi :=  \sum_{k=1}^5 x_k^5 - 5\psi\prod_{k=1}^5 x_k  .
\]
Let \mufive\ be the multiplicative group of $5^{\text{th}}$ roots
of unity, and let
\[
\widetilde{G}:=\{\vec{\alpha}=(\alpha_1,\dots,\alpha_5)\in
(\mufive)^5 \suchthat
\prod_{k=1}^5 \alpha_k = 1\}
\]
act  on ${\Bbb P}^4$ by $\vec{\alpha}: x_i\mapsto \alpha_i\cdot x_i$.
There is a ``diagonal'' subgroup of order $5$ which acts trivially; let
$G=\widetilde{G}/\{\mbox{diagonal}\}$
be the image of $\widetilde{G}$ in $\Aut({\Bbb P}^4)$.  $G$ is a
group which is abstractly isomorphic to $({\Bbb Z}/5{\Bbb Z})^3$.

The action of $G$ preserves the threefold ${\cal Q}_\psi$; let
$\eta:{\cal Q}_\psi\to {\cal Q}_\psi/G$ denote the quotient map.
For each pair of distinct indices $i$, $j$, the set of 5 points
\[
S_{ij}:=\{x_i^5+x_j^5=0, x_\ell =0 \text{ for all }
 \ell\ne i,j\}\subset {\cal Q}_\psi
\]
 is preserved by $G$, and there
is a group $G_{ij}\subset G$ of order 25 which is the stabilizer
of each point in the set.  The image $S_{ij}/G$
is a single point $p_{ij}\in {\cal Q}_\psi/G$.  In addition, for each triple
of distinct indices $i$, $j$, $k$, the curve
\[
\widetilde C_{ijk}:=\{x_i^5+x_j^5+x_k^5=0, x_\ell =0  \text{ for all }
 \ell \ne i,j,k \}\subset {\cal Q}_\psi
\]
is preserved by $G$.  There is a subgroup $G_{ijk}\subset G$ of order
5 which is the stabilizer of every point in
$\widetilde C - \eta^{-1}(\{p_{ij},p_{jk},p_{ik}\})$.
The image $C_{ijk}=\widetilde C_{ijk}/G$ is a smooth curve in
${\cal Q}_\psi/G$.
The action of $G$ is free away from the curves $\widetilde C_{ijk}$.

The quotient
space ${\cal Q}_\psi/G$ has only canonical singularities.
At most points of $C_{ijk}$, the surface section of the singularity is a
rational
double point of type $A_4$, but at the points $p_{ij}$ the singularity
is more complicated:
three of the curves of $A_4$-singularities meet at each $p_{ij}$.
By a theorem of Markushevich \cite[Prop.~4]{markushevich} and
Roan \cite[Prop.~2]{roan1}, these singularities can be
resolved to give a Calabi-Yau manifold ${\cal W}_\psi$.
There are choices to be made in this resolution
process; we describe a particular
choice in appendix B.
(By another theorem of Roan \cite[Lemma 4]{roan2},
 any two resolutions differ by a sequence of flops.)

For any $\alpha\in\mufive$,
there is a natural isomorphism between ${\cal Q}_{\alpha\psi}/G$ and
${\cal Q}_{\psi}/G$ induced by the map
\begin{equation}
\label{isomorphism}
(x_1,x_2,x_3,x_4,x_5)\mapsto(\alpha^{-1}x_1,x_2,x_3,x_4,x_5).
\end{equation}
This extends to an isomorphism between ${\cal W}_{\alpha\psi}$ and ${\cal
W}_\psi$,
provided that we have resolved singularities in a compatible way.
We verify in appendix B that the choices in the resolution can
be made in a sufficiently natural way that this isomorphism is
guaranteed to exist.

Thus, $\lambda:=\psi^5$ is a more natural parameter to use for our family.
We define the {\em quintic-mirror family} to be
\[
\{ {\cal W}_{\sqrt[5]{\lambda}} \} \to \{\lambda\}
\cong {\Bbb C}.
\]
This has a natural compactification to a family over ${\Bbb P}^1$,
with boundary $B={\Bbb P}^1-{\Bbb C}=\{\infty\}$.

The computation made by Candelas et al.\ \cite{pair} shows that
the monodromy at
$\infty$ is maximally unipotent,
and that $m=1$, i.e., that the monodromy is small.
(We explain in appendix C how this follows from \cite{pair}.)
The key computation in \cite{pair} is an explicit calculation of
the $q$-expansion of the mathematically normalized Yukawa coupling.
Candelas et al.\ find that the $q$-expansion begins:
\begin{equation} \label{formula1}
\kappa_{ttt} = 5 + 2875e^{2\pi it} + 4876875e^{4\pi it} + \cdots.
\end{equation}
In fact, they have computed at least 10 coefficients.

\section{Mirror moonshine?}
\label{moonshine}

Greene and Plesser \cite{greene-plesser}, using arguments from
superconformal field
 theory, have identified the
family of quintic-mirrors $\{ {\cal W}_{\sqrt[5]{\lambda}} \}$
as the ``mirror'' of the  family of smooth
quintic threefolds $\{ {\cal M}_z \}$.  Note that
 the Hodge numbers satisfy
\[\begin{array}{ll}
h^{1,1}({\cal M})=1 &
h^{2,1}({\cal M})=101\\
h^{1,1}({\cal W})=101 &
h^{2,1}({\cal W})=1.
\end{array}\]
According to the mirror symmetry conjecture,
 varying the complex structure in the family $\{ {\cal W}_{\sqrt[5]{\lambda}}
\}$
 should
correspond to varying the ``extra structure'' $S$ on a fixed smooth quintic
threefold ${\cal M}$.  These are both one-parameter variations.

Candelas et al.\ \cite{pair}, arguing from physical principles,
 propose an identification of the Yukawa coupling
 of the quintic-mirrors
with the cubic form on $H^{1,1}({\cal M})$.
In terms of the mathematical framework established here, that identification
involves four assertions:
\begin{enumerate}
\item[(i)]
The isomorphism $H^{-1,1}({\cal W})\to H^{1,1}({\cal M})$  defined by
$d/dt\mapsto [H]$
(where $d/dt\in H^{-1,1}({\cal W})$ is the vector
field defined by the canonical parameter $t$, and $[H]\in H^{1,1}({\cal M})$
 is the class
of a hyperplane section of ${\cal M}$) is the isomorphism
which is
 predicted by the mathematical version of the mirror symmetry
conjecture.

\item[(ii)]
The mathematically normalized Yukawa coupling $\kappa_{ttt}$ on
$H^{-1,1}({\cal W})$ is
the correctly normalized coupling predicted by the physical theory.

\item[(iii)]
The parameter $R$ from the physical theory
coincides with $-2\pi it$, where $t$ is again the canonical
parameter.
Thus, the $q$-expansion of $\kappa_{ttt}$ in equation (\ref{qexpt}) will
coincide with the asymptotic expansion in $R$ given by equation (\ref{asymp}),
evaluated on the generator $H^{\otimes 3}$ of $\Sym^3H^{-1,1}({\cal M})$.

\item[(iv)]
There is an explicit formula for the coefficients $\sigma_k$, as described
below.

\end{enumerate}

To explain the formula for $\sigma_k$, let $n_k$ denote
 the number of rational
curves of degree $k$ on the generic quintic threefold.
Candelas et al.\ propose the formula
\begin{equation}  \label{formula2}
\kappa_{ttt} = 5 +
\sum_{k=1}^\infty \frac{n_kk^3e^{2\pi ikt}}{1-e^{2\pi ikt}}
= 5 + n_1e^{2\pi it} + (2^3n_2 + n_1)e^{4\pi it} + \cdots,
\end{equation}
which implicitly incorporates their
expressions
for the higher coefficients
(The first two expressions are $\sigma_1(H^{\otimes 3})=n_1$,
 $\sigma_2(H^{\otimes 3})=2^3n_2 + n_1$.)

In the large radius limit $\Im t \to\infty$, the right hand side of
equation (\ref{formula2}) approaches $5$.  This agrees with the
mirror symmetry conjecture,\footnote{This should not be taken as strong
evidence in favor of the conjecture, since the definitions have been
carefully designed to ensure that this limit would be correct.}
since the topological intersection form
on ${\cal M}$ is determined by its value on the standard generator $H$,
viz., $H^3=5$.

Moreover, by
comparing equations (\ref{formula1}) and (\ref{formula2}), we
can {\em predict\/} values for the numbers $n_k$.  The first two
predictions are
$n_1=2875$, which was classically known
to be the number of lines on a quintic threefold, and $n_2=609250$,
which coincides with  the number of conics on a quintic
threefold
computed by Katz \cite{katz}!

\raisebox{1.2ex}{\makebox[0pt][l]{\underline{\phantom
  {Unfortunately, there seem to be difficulties with $n_3$.
  }}}}
Unfortunately, there seem to be difficulties with $n_3$.
\quad Not any more!!

\bigskip

How was formula (\ref{formula2}) arrived at?  I am told that the field
theory computation necessary to derive this formula can be
done in principle, but seems to be too hard to carry out in practice at
present.  So Candelas et al.\ give a rough derivation of this formula
based on some assumptions.
Why do they believe the resulting formula to be correct?  I quote from
\cite{pair}:
\begin{quote}
These numbers provide compelling evidence that our assumption about the
form of the prefactor is in fact correct.  The evidence is not so much
that we obtain in this way the correct values for $n_1$ and $n_2$,
but rather that the coefficients in eq.~(\ref{formula1}) have remarkable
divisibility properties.  For example asserting that the second coefficient
$4,876,875$ is of the form $2^3n_2+n_1$ requires that the result of
subtracting $n_1$ from the coefficient yields an integer that is divisible
by $2^3$.  Similarly, the result of subtracting $n_1$ from the third
coefficient must yield an integer divisible by $3^3$.  These conditions
become increasingly intricate for large $k$.  It is therefore remarkable
that the $n_k$ calculated in this way turn out to be integers.
\end{quote}
I would add that it is equally remarkable that the coefficients in
eq.~(\ref{formula1}) themselves turn out to be integers:  I know of
no proof that this is the case.

These arguments  have
 a rather numerological flavor.
I am reminded of the
 numerological observations made by Thompson
\cite{numerology} and Conway and Norton \cite{monster}
about the $j$-function and the monster group.
At the time those papers were written, no connection between these
two mathematical objects was known.  The $q$-expansion of the
$j$-function was known to have integer coefficients, and it was
observed that these integers were integral linear combinations
of the degrees of irreducible representations of the monster group.
This prompted much speculation about possible deep connections
between the two, but at the outset all such speculation had to
be characterized as ``moonshine'' (Conway and Norton's term).

The formal similarities to the present work should be clear:
a $q$-expansion
of some kind is found to have integer coefficients, and these integers
then appear to be integral linear combinations of another set of integers,
which
occur elsewhere in mathematics in a rather unexpected location.
Perhaps it is too much to hope that the eventual explanation will be
as pretty in this case.

\section*{Appendix A:  Proofs of the monodromy lemmas}

Let
\[
W_0 \subset W_1 \subset \dots \subset W_{2n}=H^n(X_{P'},{\Bbb Q})
\]
be the monodromy weight filtration at $P$, and let
\[
 F^0 \supset F^1 \supset \cdots \supset
F^{n-1} \supset
F^n \supset  (0).
\]
be the limiting Hodge filtration at $P$.
(We refer the reader to \cite{transcendental} or \cite{schmid} for the
definitions.)
By a theorem of Schmid
\cite{schmid}, these induce a mixed Hodge structure on the cohomology.
Note that since $N^{n+1}=0$, we have $W_0=\Im N^n$.

Moreover, if $\ip{\ }{\ }$ denotes
the polarization on the cohomology, we have
\[
\ip{Nx}{y}=-\ip{x}{Ny}.
\]
Recall also that the polarization is symmetric or skew-symmetric,
depending on the dimension $n$:
\[
\ip{x}{y}=(-1)^n\ip{y}{x}
\]

\begin{pf*}{Proof of lemma \ref{lemma1}}
Since $W_{\textstyle\cdot}$
is the monodromy weight filtration, $N^n$ induces
 an isomorphism
\begin{equation}  \label{eq1}
N^n: W_{2n}/W_{2n-1} \to W_0.
\end{equation}
These spaces cannot be zero, since $(T_P-I)^n\ne 0$.  On the other hand,
since $F^{n+1}=(0)$, the Hodge structure on $W_{2n}/W_{2n-1}$ must be
purely of type $(n,n)$.  It follows that
$F^n/(F^n\cap W_{2n-1})=W_{2n}/W_{2n-1}$.  But since $F^n$ is one-dimensional,
this can only happen if $F^n \subset W_{2n}-W_{2n-1}$, and
$W_{2n}/W_{2n-1}$ has dimension one.  By the isomorphism (\ref{eq1}),
$W_0=\Im N^n$ has dimension one as well.

Next, note that $W_{2n-1}/W_{2n-2}$ has a Hodge structure with two types,
$(n,n-1)$ and $(n-1,n)$, each of which must determine a space of half
the total dimension.  But since $F^n\cap W_{2n-1} = (0)$, nothing non-zero
can have type $(n,n-1)$.  It follows that $W_{2n-1}/W_{2n-2}=(0)$,
and that $W_1/W_0=(0)$ as well (using the isomorphism induced by
$N^{n-1}$).

Thus, the image of $N^{n-1}$ comes entirely from the map
\[ N^{n-1}:W_{2n}\to W_2.\]
That this image is two-dimensional is easily seen:
$W_0$ is one-dimensional, and there is an isomorphism
\[ N^{n-1}:W_{2n}/W_{2n-1}\to (\Im N^{n-1})/W_0, \]
which shows that $(\Im N^{n-1})/W_0$ is also one-dimensional.
\end{pf*}

In order to prove lemma \ref{lemma2}, we must first prove

\begin{lemma}[Essentially due to Friedman and Scattone \cite{fs}]
\label{lemma3}
\mbox{}

\noindent
Good integral bases exist, and form bases of the two-dimensional
${\Bbb Q}$-vector space $\Im N^{n-1}$.
If $g_0$, $g_1=\frac1\lambda N^{n-1}g$ is a good integral basis, then
\begin{equation}
\label{nn-1}
\frac1\lambda N^{n-1}x= - \ip{g_1}{x}g_0 + \ip{g_0}{x}g_1
\end{equation}
for all $x\in H^n(X_{P'},{\Bbb Q})$.
\end{lemma}

\begin{pf}  Choose either generator of $\Im N^n \cap H^n(X_{P'},{\Bbb Z})$
as $g_0$.
We claim that $(g_0)^\perp = W_{2n-2}$.  Let $h\in W_{2n}-W_{2n-2}$,
so that $N^nh=ag_0$ with $a\ne 0$.  Then for any $x$ we have
\[
\ip{N^nx}{h}=(-1)^n\ip{x}{N^nh}=(-1)^na\ip{x}{g_0}.
\]
Thus, $W_{2n-2}=\ker N^n \subset (g_0)^\perp$.
Since both $W_{2n-2}$ and $(g_0)^\perp$ are codimension one subspaces
of $W_{2n}$, they must be equal.

By Poincar\'e duality, the polarization on
$H^n(X_{P'},{\Bbb Z})$ is a unimodular pairing.
Thus, there exists an element $g\in H^n(X_{P'},{\Bbb Z})$ such that
$\ip{g_0}{g}=1$.
Since $g\not\in(g_0)^\perp$, neither $N^{n-1}g$ nor $N^n g$ is zero.
There is thus a unique positive rational number $\lambda$
such that $g_1 = \frac1\lambda N^{n-1}g$ is an indivisible element of
$H^n(X_{P'},{\Bbb Z})$.
It is clear that $g_0$, $g_1$ forms
a basis for the ${\Bbb Q}$-vector space $\Im N^{n-1}$.

We next claim that $\ip{g_1}{g}=\frac1\lambda \ip{N^{n-1}g}{g}=0$.
For on the one hand, moving the $N$'s
to the right side one at a time we have
\[
\ip{N^{n-1}g}{g}=(-1)^{n-1}\ip{g}{N^{n-1}g}
\]
while on the other hand, the symmetry of the polarization says that
\[
\ip{N^{n-1}g}{g}=(-1)^{n}\ip{g}{N^{n-1}g}.
\]
It follows that $\ip{N^{n-1}g}{g}=0$.

To prove equation (\ref{nn-1}), we first compute in general
\[
\ip{N^{n-1}x}{g}=(-1)^{n-1}\ip{x}{N^{n-1}g}=(-1)^{n-1}\ip{x}{\lambda g_1}
=-\lambda\ip{g_1}{x}.
\]
Now suppose that $x\in W_{2n-2}$.
Then $N^{n-1}x\in \Im N^n$, which implies that $N^{n-1}x=ag_0$ for some
$a$.  Thus, in this case
\[
\ip{N^{n-1}x}{g} = \ip{ag_0}{g} = a,
\]
which implies that $a=-\lambda\ip{g_1}{x}$.  Thus,
\[
\frac1\lambda N^{n-1}x=\frac1\lambda ag_0=\ip{g_1}{x}g_0
\]
and since $\ip{x}{g_0}=0$, the formula follows in this case.

To prove the formula in general, note that
\[
\ip{g_0}{x-\ip{g_0}{x}g\vphantom{N^{n-1}}}= 0
\]
for any $x$, so that $x-\ip{g_0}{x}g\in (g_0)^\perp=W_{2n-2}$.
Thus, applying the
previous case we find
\begin{eqnarray*}
\frac1\lambda N^{n-1}x &=& \frac1\lambda N^{n-1}(x-\ip{g_0}{x}g)
           + \frac1\lambda N^{n-1}(\ip{g_0}{x}g)\\
 &=& -\ip{g_1}{(x-\ip{g_0}{x}g)\vphantom{N^{n-1}}}g_0
           + \ip{g_0}{x}\frac1\lambda N^{n-1}g\\
 &=& -\left(\ip{g_1}{x}-\ip{g_0}{x}\ip{g_1}{g}\vphantom{N^{n-1}}\right)g_0
           + \ip{g_0}{x} g_1\\
 &=& -\ip{g_1}{x}g_0 + \ip{g_0}{x}g_1
\end{eqnarray*}
since $\ip{g_1}{g}=0$.
\end{pf}

We can now prove lemma \ref{lemma2}.

\begin{pf*}{Proof of lemma \ref{lemma2}}
The only generators of $\Im N^n \cap H^n(X_{P'},{\Bbb Z})$ are $\pm g_0$,
so we must have $g_0'=(-1)^\ell g_0$ for some $\ell\in\Bbb Z$.
Write $g_1'=\frac1{\lambda'}N^{n-1}g'$ for some $g'$ with
$\ip{g_0'}{g'}=1$, and let $k=-\ip{g_1}{g'}\in\Bbb Z$.  Then
by lemma \ref{lemma3},
\[
\frac1\lambda N^{n-1}g'  =  -\ip{g_1}{g'}g_0+\ip{g_0}{g'}g_1
  =  k\,g_0+(-1)^\ell g_1.
\]
Thus, $\frac1\lambda N^{n-1}g' \in H^n(X_{P'},{\Bbb Z})$.
We claim that it must be an indivisible element there.  For if
$\frac1{\lambda\mu} N^{n-1}g'$ is integral for some $\mu\in\Bbb Z$
with $\mu>1$,
then reversing the roles of $g$ and $g'$ in the argument above
shows that $\frac1{\lambda\mu} N^{n-1} g$ is also integral, a
contradiction.

Thus, $g_1'=k\,g_0+(-1)^\ell g_1$.
\end{pf*}

\section*{Appendix B:  Resolutions of certain quotient singularities}

In this appendix, we will verify that the singularities of the variety
${\cal Q}_\psi/G$ can be resolved in a natural way.  The choices we make are
sufficiently natural that the isomorphism between ${\cal Q}_{\alpha\psi}/G$
and ${\cal Q}_\psi/G$ automatically lifts to an isomorphism between the
desingularizations.

We choose to follow the strategy outlined by Reid \cite{reid} for
resolving canonical threefold singularities.  In brief, we perform
the following steps:
\begin{list}{}{\advance\leftmargin by 1.6em \advance\labelwidth by 1.6em}
\item[Step I:]
Blow up the ``non-cDV points'' of ${\cal Q}_\psi/G$.  (These are exactly
the 10 points $p_{ij}\in {\cal Q}_\psi/G$ which are the images of
points in ${\cal Q}_\psi$ with stabilizer of order 25.)

\item[Step IIA:]
Blow up the singular locus.  (It has pure dimension one.)

\item[Step IIB:]
Blow up the pure dimension one part of the singular locus.  (60 isolated
singular points (lying over the $p_{ij}$)
were created by step IIA, and these are not to be
blown up yet.)

\item[Step III:]
Obtain a projective small resolution of the remaining 60 singular points
by blowing up the union of the proper transforms of the exceptional
divisors from step I.
\end{list}
Step III involved an additional choice, since Reid's strategy does not
specify how one should obtain small resolutions.

When stated in this form, it is clear that the process we have described
is sufficiently natural that it is preserved under any isomorphism.
It yields a projective (hence K\"ahler) variety ${\cal W}_\psi$ with trivial
canonical bundle.

In the remainder of this appendix, we will show that the process above
has the properties mentioned during its description, and that
it gives a resolution of singularities of ${\cal Q}_\psi/G$.

We first observe the effect of the process on the curve $C_{ijk}$,
away from the points $p_{ij}$, $p_{jk}$, $p_{ik}$.
Steps I and III are concentrated at
those special points (and their inverse images) and so these steps do
not affect $C_{ijk}$.  Steps IIA and IIB simply blow up $C_{ijk}$ and
then the residual singular curve in the exceptional divisor.
But two blowups are precisely what is
required to resolve a rational double point of type $A_4$, as is
easily verified from its equation $xy+z^5=0$.

To verify that the process has the correct properties at the points $p_{ij}$,
we use the language of toroidal embeddings (see \cite{te} or \cite{oda}).
It suffices to consider the point $p_{45}$.  Since $(x_1,x_2,x_3)$
serve as coordinates in a neighborhood of any of the points in
$\eta^{-1}(p_{45})$, the
singularity $p_{45}\in {\cal Q}_\psi/G$ is isomorphic to a neighborhood of
the origin in $\Bbb C^3/G_{45}$, where
$G_{45}\cong\{(\alpha_1,\alpha_2,\alpha_3)
\in(\mufive)^3\suchthat\prod\alpha_k=1\}$
acts diagonally on $\Bbb C^3$.

Let $M$ be the lattice of $G_{45}$-invariant rational
monomials in $\Bbb C(x_1,x_2,x_3)$.
We embed $M$ in $\Bbb R^3$ by\footnote{This non-standard embedding is
chosen in order to make the coordinates of the dual lattice be integers.}
\[
M=\{(m_1,m_2,m_3)\in\Bbb R^3\suchthat
x_1^{5m_1}x_2^{5m_2}x_3^{5m_3}\in\Bbb C(x_1,x_2,x_3)^{G_{45}}\}.
\]
It is easy to see that $\{(1,0,0),(0,1,0),(\frac15,\frac15,\frac15)\}$
is a basis of the lattice $M\subset\Bbb R^3$.
Let
\begin{eqnarray*}
N & = & \{\vec{n}\in\Bbb R^3\suchthat \vec{m}\cdot\vec{n}\in\Bbb Z
\text{\ for all\ }
\vec{m}\in M\} \\
 & = & \{(n_1,n_2,n_3)\in\Bbb Z^3\suchthat n_1+n_2+n_3\equiv 0 \mod 5\}
\end{eqnarray*}
be the dual lattice, and let $\sigma\subset N_{\Bbb R}$ be the convex
cone generated by $(5,0,0)$, $(0,5,0)$, and $(0,0,5)$.
According to the theory of toroidal embeddings,
\[
\Bbb C^3/G_{45} = \Spec\Bbb C[x_1,x_2,x_3]^{G_{45}} = U_\sigma,
\]
where $U_\sigma$ is the toric variety associated to $\sigma$.

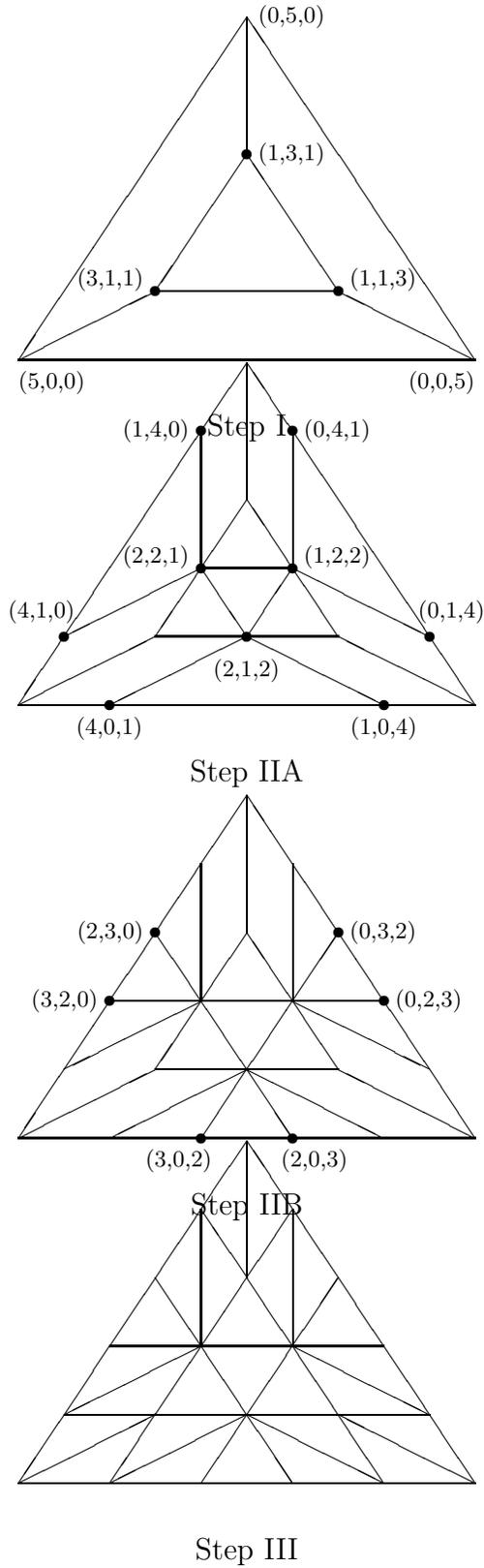
\begin{figure}[t]
\setlength{\unitlength}{.38em}
{
\bigskip

\noindent
\begin{picture}(40,30)

\put(20,-6){\makebox(0,0){Step I}}

\thinlines

\put(0,0){\line(1,0){40}}
\put(0,-3){\makebox(0,2)[l]{\scriptsize (5,0,0)}}

\put(40,0){\line(-2,3){20}}
\put(40,-3){\makebox(0,2)[r]{\scriptsize (0,0,5)}}

\put(20,30){\line(-2,-3){20}}
\put(21,30){\makebox(0,0)[l]{\scriptsize (0,5,0)}}

\put(12,6){\line(1,0){16}}
\put(12,6){\circle*{1}}
\put(11,6){\makebox(0,2)[r]{\scriptsize (3,1,1)}}

\put(28,6){\line(-2,3){8}}
\put(28,6){\circle*{1}}
\put(29,6){\makebox(0,2)[l]{\scriptsize (1,1,3)}}

\put(20,18){\line(-2,-3){8}}
\put(20,18){\circle*{1}}
\put(21,18){\makebox(0,0)[l]{\scriptsize (1,3,1)}}

\put(0,0){\line(2,1){12}}
\put(40,0){\line(-2,1){12}}
\put(20,30){\line(0,-1){12}}

\end{picture}
\hfill
\begin{picture}(40,30)

\put(20,-6){\makebox(0,0){Step IIA}}

\thinlines

\put(0,0){\line(1,0){40}}
\put(40,0){\line(-2,3){20}}
\put(20,30){\line(-2,-3){20}}

\put(12,6){\line(1,0){16}}
\put(28,6){\line(-2,3){8}}
\put(20,18){\line(-2,-3){8}}

\put(0,0){\line(2,1){12}}
\put(40,0){\line(-2,1){12}}
\put(20,30){\line(0,-1){12}}

\put(20,6){\line(2,3){4}}
\put(20,6){\circle*{1}}
\put(20,2){\makebox(0,2){\scriptsize (2,1,2)}}
\put(20,6){\line(-2,-1){12}}
\put(20,6){\line(2,-1){12}}

\put(24,12){\line(-1,0){8}}
\put(24,12){\circle*{1}}
\put(25,12){\makebox(0,2)[l]{\scriptsize (1,2,2)}}
\put(24,12){\line(2,-1){12}}
\put(24,12){\line(0,1){12}}

\put(16,12){\line(2,-3){4}}
\put(16,12){\circle*{1}}
\put(15,12){\makebox(0,2)[r]{\scriptsize (2,2,1)}}
\put(16,12){\line(0,1){12}}
\put(16,12){\line(-2,-1){12}}

\put(8,0){\circle*{1}}
\put(8,-3){\makebox(0,2){\scriptsize (4,0,1)}}

\put(32,0){\circle*{1}}
\put(32,-3){\makebox(0,2){\scriptsize (1,0,4)}}

\put(36,6){\circle*{1}}
\put(35,8.2){\makebox(0,0)[l]{\scriptsize (0,1,4)}}

\put(24,24){\circle*{1}}
\put(25,24){\makebox(0,0)[l]{\scriptsize (0,4,1)}}

\put(16,24){\circle*{1}}
\put(15,24){\makebox(0,0)[r]{\scriptsize (1,4,0)}}

\put(4,6){\circle*{1}}
\put(5,8.2){\makebox(0,0)[r]{\scriptsize (4,1,0)}}

\end{picture}

\bigskip

\bigskip

\bigskip

\medskip

\noindent
\begin{picture}(40,30)

\put(20,-6){\makebox(0,0){Step IIB}}

\thinlines

\put(0,0){\line(1,0){40}}
\put(40,0){\line(-2,3){20}}
\put(20,30){\line(-2,-3){20}}

\put(12,6){\line(1,0){16}}
\put(28,6){\line(-2,3){8}}
\put(20,18){\line(-2,-3){8}}

\put(0,0){\line(2,1){12}}
\put(40,0){\line(-2,1){12}}
\put(20,30){\line(0,-1){12}}

\put(16,0){\line(2,3){12}}

\put(20,6){\line(-2,-1){12}}
\put(20,6){\line(2,-1){12}}

\put(32,12){\line(-1,0){24}}

\put(24,12){\line(2,-1){12}}
\put(24,12){\line(0,1){12}}

\put(12,18){\line(2,-3){12}}

\put(16,12){\line(0,1){12}}
\put(16,12){\line(-2,-1){12}}

\put(16,0){\circle*{1}}
\put(17,-3){\makebox(0,2)[r]{\scriptsize (3,0,2)}}

\put(24,0){\circle*{1}}
\put(23,-3){\makebox(0,2)[l]{\scriptsize (2,0,3)}}

\put(32,12){\circle*{1}}
\put(33,12){\makebox(0,0)[l]{\scriptsize (0,2,3)}}
\put(28,18){\circle*{1}}
\put(29,18){\makebox(0,0)[l]{\scriptsize (0,3,2)}}

\put(12,18){\circle*{1}}
\put(11,18){\makebox(0,0)[r]{\scriptsize (2,3,0)}}
\put(8,12){\circle*{1}}
\put(7,12){\makebox(0,0)[r]{\scriptsize (3,2,0)}}

\end{picture}
\hfill
\begin{picture}(40,30)

\put(20,-6){\makebox(0,0){Step III}}

\thinlines

\put(0,0){\line(1,0){40}}
\put(40,0){\line(-2,3){20}}
\put(20,30){\line(-2,-3){20}}

\put(4,6){\line(1,0){32}}
\put(32,0){\line(-2,3){16}}
\put(24,24){\line(-2,-3){16}}

\put(0,0){\line(2,1){12}}
\put(40,0){\line(-2,1){12}}
\put(20,30){\line(0,-1){12}}

\put(16,0){\line(2,3){12}}

\put(20,6){\line(-2,-1){12}}
\put(20,6){\line(2,-1){12}}

\put(32,12){\line(-1,0){24}}

\put(24,12){\line(2,-1){12}}
\put(24,12){\line(0,1){12}}

\put(12,18){\line(2,-3){12}}

\put(16,12){\line(0,1){12}}
\put(16,12){\line(-2,-1){12}}

\end{picture}

\bigskip

\bigskip
}
\caption{The steps in the toroidal resolution. \protect\rule[-3ex]{0pt}{3ex}}
\label{fig1}
\end{figure}

Each blowup of $U_\sigma$ corresponds to a decomposition of $\sigma$
into a fan.  The effects of the blowups in our process is illustrated
in figure \ref{fig1}, which depicts the intersection of the fan with
$\{(n_1,n_2,n_3)\in N_{\Bbb R}\suchthat\sum n_k=5\}$ after each step.
The exceptional divisors $D_{\vec{n}}$ of each blowup are indicated
by solid dots, labeled by the corresponding elements $\vec{n}\in N$.
(The fact that the stated blowups produce the illustrated decomposition
is a straightforward calculation with the toroidal embeddings.)

We can now see in detail what happens in our process.  In
step I, we blow up $p_{45}$, and produce three exceptional divisors
$D_{(3,1,1)}$, $D_{(1,3,1)}$, and $D_{(1,1,3)}$.
The remaining singular locus at this stage consists of the original
three curves of $A_4$-singularities together with three new curves
of $A_1$-singularities:  the intersections of pairs of exceptional
divisors.  In step IIA, we blow up the union of these six curves, and
produce nine new exceptional divisors:  one corresponding to each
curve of $A_1$-singularities (such as $D_{(2,1,2)}$), and two
corresponding to each curve of $A_4$-singularities (such as $D_{(4,0,1)}$
and $D_{(1,0,4)}$).  The remaining singularities consist of six
isolated points (corresponding to the quadrilaterals in the figure)
and three curves:  the intersections of the corresponding pairs of
exceptional divisors from the original $A_4$-singularities.

In step IIB, we blow up these three curves, producing six new
exceptional divisors, two for each curve
(such as $D_{(3,0,2)}$ and $D_{(2,0,3)}$).
This leaves the six isolated singular points; but blowing up the
proper transforms of $D_{(3,1,1)}$, $D_{(1,3,1)}$, and $D_{(1,1,3)}$
(which are now disjoint) in step III resolves those final singular points.

\section*{Appendix C:  The monodromy of the quintic-mirrors}

In this appendix we will explain how to use the calculation of
Candelas et al.\ \cite{pair} to verify the monodromy statements
about the family of quintic-mirrors which we made in section \ref{qm}.

Candelas et al.\ begin by choosing an explicit
basis $\{A^1,A^2,B_1,B_2\}$
for the homology $H_3({\cal W}_\psi,\Bbb Z)$ of a quintic-mirror,
valid in some simply-connected region in $\{\psi\suchthat\psi^5\ne 0,1\}$
which includes the wedge
$\{\psi\suchthat 0<\arg\psi <2\pi/5\}$.
This basis is {\em symplectic}, i.e., ${A^a}\cap{B_b}={\delta^a}_b$
and ${A^a}\cap{A^b}={B_a}\cap{B_b}=0$.
The corresponding dual basis of $H^3({\cal W}_\psi,\Bbb Z)$
is denoted by $\{\alpha_1,\alpha_2,\beta^1,\beta^2\}$.
Fixing a particular holomorphic 3-form
$\Omega$ (which depends on $\psi$), we then get {\em period functions\/}
\[
z^a=\int_{A^a}\Omega ,\ \ \ \ {\cal G}_b=\int_{B_b}\Omega.
\]
These fit into a {\em period vector\/}
\[
\Pi = \left( \begin{array}{c} {\cal G}_1 \\ {\cal G}_2 \\  z^1 \\ z^2
\end{array} \right) .
\]
By doing some integrals, calculating the differential equation
satisfied by a period function, and manipulating certain hypergeometric
functions, the authors of \cite{pair} are able to obtain explicit formulas
for the four period functions.  This allows them to calculate
the monodromy of the periods around various paths.

Notice that we are working in the $\psi$-plane at present.  The family
$\{{\cal W}_\psi\}$ has singular fibers at $\psi=0$ and at $\psi=\alpha$
for all fifth roots of unity $\alpha$; there is also a singular fiber
over $\psi=\infty$.  Candelas et al.\ calculate the monodromy on the
periods induced by transport around $\psi=1$, which they represent in
matrix form by $\Pi\to T\Pi$.  They also compute, for $|\psi|<1$,
 the effect on the
periods of the isomorphism ${\cal W}_{\alpha\psi} \cong {\cal W}_\psi$ given in
equation (\ref{isomorphism}), representing this by
$\Pi(\alpha\psi)=A\Pi(\psi)$.

We need to know the monodromy around $\infty$ in the $\lambda$-plane,
where $\lambda=\psi^5$.  A moment's thought will convince the reader
that this is represented by
\[
\Pi\to (T^{-1}A^{-1})\Pi,
\]
and that $(AT)^{-5}$ describes the monodromy around $\infty$ in the
$\psi$-plane (as asserted in \cite{pair}).
Let $T_P=T^{-1}A^{-1}$.

The explicit calculations from \cite{pair} for the matrices $A$ and $T$ are:
\[
A=
\left( \begin{array}{rrrr}
-9 & -3 & 5 & 3 \\
0 & 1 & 0 & -1 \\
-20 & -5 & 11 & 5 \\
-15 & 5 & 8 & -4
\end{array} \right)
\ \ \ \
T=
\left( \begin{array}{rrrr}
1 & 0 & 0 & 0 \\
0 & 1 & 0 & 1 \\
0 & 0 & 1 & 0 \\
0 & 0 & 0 & 1
\end{array} \right)
\]
from which it easily follows that
\[
(\log (T_P))^2=
\left( \begin{array}{rrrr}
0 & 5 & 0 & 0 \\
0 & 0 & 0 & 0 \\
0 & 10 & 0 & 0 \\
-10 & 0 & 5 & 0
\end{array} \right) ,
\ \ \ \
(\log (T_P))^3=
\left( \begin{array}{rrrr}
0 & 0 & 0 & 0 \\
0 & 0 & 0 & 0 \\
0 & 0 & 0 & 0 \\
0 & -5 & 0 & 0
\end{array} \right) .
\]
In particular, the index of nilpotency
of $\log (T_P)$
is maximal.

(We note in passing that at
$\lambda=1$ the monodromy is represented by $T$, and since
$(T-I)^2=0$, the index is not maximal there.  In addition, at
$\lambda=0$ the monodromy is represented by $A$.  This monodromy
matrix is only quasi-unipotent, with $A^5=I$ unipotent; the index of $A^5$
is not maximal either.  It follows that $\lambda=\infty$ is the only
possible boundary point with maximally unipotent monodromy.)

In order to construct a good integral basis $g_0$, $g_1$, we compute
\begin{eqnarray*}
(\log (T_P))^2(\int_{A^2}\Omega) & = &
-10\int_{B_1}\Omega + 5\int_{A^1}\Omega \\
(\log (T_P))^3(\int_{A^2}\Omega) & = &
-5\int_{B_2}\Omega .
\end{eqnarray*}
Using the relations
\begin{equation}
\label{relations}
\ad_{\ip{\,}{\,}}(\alpha_a)=\int_{B_a}\ ,\ \ \ \
\ad_{\ip{\,}{\,}}(\beta^b)=-\int_{A^b}\ ,
\end{equation}
this implies
\begin{eqnarray*}
(\log (T_P))^2(\beta^2) & = &
10\alpha_1+5\beta^1 \\
(\log (T_P))^3(\beta^2) & = &
5\alpha_2 .
\end{eqnarray*}
Thus, we may take $g_0=\alpha_2$.  If we then choose $g=\beta^2$
so that $\ip{g_0}{g}=\ip{\alpha_2}{\beta^2}=1$,  we get
$\lambda=5$ and $g_1=2\alpha_1+\beta^1$.  It follows that
\[
(\log (T_P))(g_1) = \frac15(\log (T_P))^3(\beta^2)
=\alpha_2=g_0,
\]
which implies that $m=1$.  Using the relations (\ref{relations}) again,
it follows that
$\gamma_0=B_2$, $\gamma_1=2B_1-A^1$.  Thus,
$t=(\int_{2B^1-A_1}\Omega)/(\int_{A^2}\Omega)$.

We need to verify that our parameter $t$ is the same one used by
Candelas et al.  Their parameter is defined in \cite[(5.9)]{pair}
by $t=w^1/w^2$, with
 $w^1$ and $w^2$  determined by a pair of
equations
\[
\amalg = N\,\Pi,\ \ \ \ w^2={\cal G}_2,
\]
where\footnote{We
have taken the liberty of correcting a typographical error in $N$
when transcribing it from \cite{pair}.}
\[
\amalg=
\left( \begin{array}{c}
{\cal F}_1 \\ {\cal F}_2 \\ w^1 \\ w^2
\end{array} \right)
\text{\ \ and\ \ }
N=
\left( \begin{array}{rrrr}
-1 & 0 & 0 & 0 \\
0 & 0 & 0 & 1 \\
2 & 0 & -1 & 0 \\
0 & 1 & 0 & 0
\end{array} \right)
\]
represent
a vector
 $\amalg$ which is a sort of mirror analogue of the period vector $\Pi$,
and
a particular integral symplectic
matrix
$N$,
respectively.
(Sadly, in the published version of \cite{pair},
the symbols $\Pi$ and $\amalg$ were identified, making section 5.2 of that
paper difficult to read.)
 It follows that
\[
t=\frac{w^1}{w^2} = \frac{2{\cal G}_1-z^1}{{\cal G}_2}
= \frac{\int_{2B^1-A_1}\Omega}{\int_{A^2}\Omega}
\]
as required.

\section*{Acknowledgements}

It is a pleasure to acknowledge helpful conversations and e-mail exchanges with
Paul Aspinwall,
Robert Bryant,
Philip Candelas,
Brian Greene,
Yujiro Kawamata,
Ronen Plesser,
Les Saper,
Chad Schoen,
and especially Sheldon Katz
as I was struggling to understand this material.
This work was partially supported by NSF grant DMS-9103827.

\ifx\undefined\bysame
\newcommand{\bysame}{\leavevmode\hbox to3em{\hrulefill}\,}
\fi

\end{document}